\documentclass[11pt]{article}
\usepackage{amsfonts}
\usepackage{mathrsfs}
\usepackage{amsmath,amssymb}
\usepackage{mathtools}

\DeclareMathOperator{\EX}{\mathbb{E}}
\usepackage[latin1]{inputenc}
\usepackage{amsmath,amssymb}
\usepackage{latexsym}
\usepackage{epstopdf}
\usepackage{geometry}
\usepackage[active]{srcltx}
\usepackage{graphicx}
\usepackage{epstopdf}
\usepackage[
bookmarks=true,         
bookmarksnumbered=true, 
colorlinks=true, pdfstartview=FitV, linkcolor=blue, citecolor=blue,
urlcolor=blue]{hyperref}

\newtheorem {theorem}{Theorem}[section]

\newtheorem {corollary}{Corollary}[section]

\def\R{{\mathbb R}}
\def\E{{{\mathbb E}\,}}

\begin{document}

\title{Cram\'er Type Moderate Deviations for Random Fields}

\maketitle
\begin{center}
\bigskip
Aleksandr Beknazaryan$^{a}$,
Hailin Sang$^{a}$
and Yimin Xiao$^{b}$

\bigskip$^{a}$ Department of Mathematics, The University of Mississippi,
University, MS 38677, USA. E-mail: abeknaza@olemiss.edu, sang@olemiss.edu

\bigskip$^{b}$ Department of Statistics and Probability, Michigan State University, East
Lansing, MI 48824, USA. E-mail:  xiao@stt.msu.edu
\end{center}

\bigskip \textbf{Abbreviated Title: }{\Large Moderate deviations for random fields}

\begin{center}
\bigskip
\textbf{Abstract}
\end{center}
We study the Cram\'er type moderate deviation for partial sums of random fields by applying
the conjugate method. The results are applicable  to the  partial sums of  linear
random fields with short or long memory and to nonparametric regression with random field errors.\\

\noindent Keywords: Cram\'er type moderate deviation, long range dependence,
nonparametric regression, spacial linear process, random field.\\

\noindent  {\textit{MSC 2010 subject classification}: 60F10, 60G60, 62E20}

\section{Introduction}
In this paper we study the Cram\'er type moderate deviations for random fields,
in particular linear random fields (often called spatial linear processes in statistics
literature) with short or long memory (short or long range dependence). The study
of moderate deviation probabilities in non-logarithmic form for independent random
variables goes back to 1920s. The first theorem in this field was published by Khinchin
(1929) who studied a particular case of the Bernoulli random variables. In his fundamental
work, Cram\'{e}r (1938) studied the estimation of the tail probability by the standard normal
distribution under the condition that the random variable has moment generating
function in a neighborhood of the origin (cf. (\ref{Cramer}) below). This condition has been
referred to as the Cram\'{e}r condition. Cram\'{e}r's work was improved by
Petrov (1954)  (see also Petrov (1975, 1995)).  Their works have stimulated a large
amount of research on moderate and large deviations; see below for a brief (and incomplete)
review on literature related to this paper. Nowadays, the area of
moderate and large deviation deviations is not only
important in probability but also plays an important role in many applied fields,
for instance, the premium calculation problem, risk management in insurance
(cf. Asmussen and Albrecher (2010)), nonparametric estimation in statistics (see,
e.g., Bahadur and Rao (1960), van der Vaart (1998), Joutard (2006, 2013)), and
in network information theory (cf. Lee et al. (2016, 2017)).

Let $X, X_{1},X_{2},\cdots $ be a sequence of independent and identically
distributed (i.i.d.) random variables with mean $0$ and variance $\sigma^2$.
Let $S_{n}=\sum_{k=1}^n X_k$ ($n\ge 1$) be the partial sums.  By the central
limit theorem,
\begin{equation*}
\lim_{n\rightarrow \infty }\sup_{x\in \mathbb{R}}\left\vert
\mathbb P(S_{n}>x\sigma\sqrt{n})-(1-\Phi (x))\right\vert =0,
\end{equation*}
where $\Phi(x)$ is the probability distribution of the standard normal
random variable. If for a suitable sequence $c_{n}$, we have
\begin{equation}\label{Q}
\lim_{n\rightarrow \infty }\sup_{0\leq x\leq c_{n}} \bigg|\frac{\mathbb P(S_{n}>
x\sigma\sqrt{n})}{1-\Phi
(x)}-1 \bigg|=0,
\end{equation}
or $\mathbb P(S_{n}>x\sigma\sqrt{n})=(1-\Phi (x))(1+o(1))$ uniformly over $x\in \lbrack 0,\,
c_{n}]$,  then Eq. (\ref{Q}) is called moderate deviation probability or normal deviation
probability for $S_n$ since it can be estimated by the standard normal distribution. We
refer to $[0, \, c_n]$ as a range for the moderate deviation. The
most famous result  of this kind is the Cram\'{e}r type moderate deviation. Under
Cram\'{e}r's condition, one has the following Cram\'{e}r's theorem (Cram\'{e}r (1938),
Petrov (1954; 1975, p.218; or 1995, p.178)): If $x\geq 0$ and $x=o(\sqrt{n})$ then
\begin{equation} \label{2}
\frac{\mathbb P(S_{n}>x\sigma\sqrt{n})}{1-\Phi (x)}=
\exp \left\{\frac{x^{3}}{\sqrt{n}}\lambda \Big(\frac{x}{\sqrt{n}}\Big)\right\} \left[1+
O\left(\frac{x+1}{\sqrt{n}}\right)\right].
\end{equation}%
Here $\lambda (z)=\sum_{k=0}^{\infty }c_{k}z^{k}$ is a
power series with coefficients depending on the cumulants of the random
variable $X$. Eq. \eqref{2} provides more precise approximation than \eqref{Q} which holds
uniformly on the range $[0,\, c_n]$ for any $c_n = o(\sqrt{n})$. The moderate deviations
under Cram\'{e}r's condition for independent non-identically distributed random variables
were obtained by Feller (1943), Petrov (1954)  and Statulevi\v{c}ius (1966). The Cram\'{e}r
type moderate deviation has also been established for the sum of independent random
variables with $p$-th moment, $p>2$.
To name a few, for example, see Rubin and Sethuraman (1965),
Nagaev (1965, 1979), Michel (1976), Slastnikov (1978), Amosova (1979),
and Frolov (2005). It should be pointed out that the ranges the moderate deviations in these
references are smaller (e.g., $c_n = O(\sqrt{\log n})$).

The Cram\'er type moderate deviations for dependent random variables
have also been studied in the literature. Ghosh (1974), Heinrich (1990)
studied the moderate deviation for $m$-dependent random variables. Ghosh
and Babu (1977), Babu and  Singh (1978a) studied moderate deviation for
mixing processes.  Grama (1997), Grama and Haeusler (2000, 2006) and Fan, Grama and Liu (2013)
investigated the large and moderate deviations for martingales.
Babu and  Singh (1978b) established moderate deviation results
for linear processes with coefficients satisfying $\sum_{i=1}^\infty i|a_i|<\infty$.
Wu and Zhao (2008)  studied moderate deviations for stationary processes
under certain conditions in terms of the physical dependence measure. But it
can be verified that the results from Wu and Zhao (2008)  can only be
applied to linear processes with short memory and their transformations. Recently
Peligrad et al. (2013) studied the exact moderate and large deviations for
short or long memory linear processes. Sang and Xiao (2018) studied
exact moderate and large deviations for linear random fields and applied the moderate
result to prove a Davis-Gut law of the iterated logarithm. Nevertheless, in the aforementioned
works, the moderate deviations are studied for dependent random variables  with
$p$-th moment, $p>2$. The exact moderate deviation for random fields under Cram\'{e}r's
condition has not been well studied.  For example, the optimal range $[0, c_n]$ and the exact
rate of convergence in (\ref{Q}) had been unknown in the random field setting.

The main objective of this paper is to establish exact moderate deviation analogous to
\eqref{2} for random fields under Cram\'{e}r's condition. Our main result is Theorem
\ref{Theorem1} below, whose proof is based on the conjugate method to change the
probability measure as in the classical case (see, e.g., Petrov (1965, 1975)). The
extension of this method to the random field setting reveals the deep relationship
between the tail probabilities and
the properties of the cumulant generating functions of the random variables such
as the analytic radius and the bounds,  for $x$ within some ranges related to the
sum of the variances and the analytic radius of the cumulant generating functions
of these random variables. Compared with the results in Sang and Xiao (2018) for
linear random fields, Theorems \ref{Theorem1} and \ref{Theorem2} in this paper
provide more precise convergence rate in the moderate deviations and explicit information
on the range $[0, \, c_n]$, which is much bigger than the range in Theorem 2.1 in Sang
and Xiao (2018).  In Section \ref{application} we show that Theorem \ref{Theorem1}
is applicable to linear random fields with short or long memory and to nonparametric
regression analysis. The results there can be applied to approximate the quantiles
and tail conditional expectations for the partial sums of  linear random fields.

In this paper we use the following notations. For two sequences $\{a_n\}$
and $\{b_n\}$ of real numbers,
$a_{n}\mathbb{\sim}b_{n}$ means $a_{n}/b_{n}\rightarrow1$ as $n\rightarrow \infty$;
$a_{n}\propto b_n$ means that $a_{n}/b_{n}\rightarrow C$ as $n\rightarrow \infty$
for some constant $C>0$; for positive sequences, the notation $a_{n}\ll b_{n}$
or $b_n\gg a_n$ means that $a_{n}/b_{n}$
is bounded. For $d,m\in\mathbb{N}$ denote $\Gamma^d_m=[-m, m]^d
\cap \mathbb{Z}^d$. Section \ref{main} gives the main results.
In Section \ref{application} we study the application of the main results
in linear random fields and nonparametric regression. All the proofs go to Section \ref{proof}.

\bigskip

\textbf{Acknowledgement} The authors are grateful to the referee and the Associate
Editor for carefully reading the paper and for insightful suggestions that significantly
improved the presentation of the paper.  The research of Hailin Sang is supported by
the Simons Foundation Grant 586789 and the College of Liberal Arts
Faculty Grants for Research and Creative Achievement at the University of Mississippi.
The research of Yimin Xiao is partially supported by NSF grants DMS-1612885 and
DMS-1607089.

\section{Main results}\label{main}
Let $\{X_{nj}, \, n \in\mathbb{N}, j\in\mathbb{Z}^d\}$ be a random field with zero means
defined on a probability space $(\Omega, {\mathcal F}, P)$. Suppose that for each $n$, the
random variables $X_{nj}, \, j \in {\mathbb Z}^d$ are independent and satisfy the following
Cram\'er condition: There is a positive constant $H_n$ such that the
cumulant generating function
\begin{equation}\label{Cramer}
L_{nj}(z)=\log \EX e^{zX_{nj}} \ \ \hbox{ of $X_{nj}$ is analytic in } D_n,
\end{equation}
where $D_n = \{ z \in \mathbb C: |z| < H_n\}$ is the disc of radius $H_n$ on the complex plane
$\mathbb C$, and  $\log$ denotes the principal value of the logarithm so that
$L_{nj}(0)=0$. This setting is 
convenient for applications to linear random fields in Section \ref{application}.

Without loss of generality we assume in this section that $\limsup\limits_{n\to \infty}  H_n<\infty$.
Within the disc $\{z \in \mathbb C: |z|< H_n\}$, $L_{nj}$ can be expanded in a convergent
power series
$$L_{nj}(z)=\sum_{k=1}^\infty \frac{\gamma_{knj}}{k!}z^k,$$
where $\gamma_{knj}$ is the cumulant of order $k$ of the random variable $X_{nj}$.
We have that  $\gamma_{1nj}=  \EX X_{nj}=0$ and   $\gamma_{2nj}=  \EX X_{nj}^2=\sigma_{nj}^2$.
By Taylor's expansion, one can verify that a sufficient condition for (\ref{Cramer})
is the following moment condition
\[
|\E X_{nj}^m| \le \frac {m!} 2 \sigma_{nj}^2 H_n^{2-m} \quad \hbox{ for all } \, m \ge 2.
\]
This condition has been used frequently in probability and statistics, see Petrov (1975, p.55),
Johnstone (1999, p.64), Picard and Tribouley (2000, p.301), Zhang and Wong (2003, p.164), among
others.

Denote
$$
S_n=\sum_{ j\in\mathbb{Z}^d}X_{nj},\quad  S_{m,n}=\sum_{j\in\Gamma^d_m} X_{nj},
$$

$$ B_n=\sum_{ j\in\mathbb{Z}^d} \sigma_{nj}^2,
\quad\mathrm{  }\quad  F_n(x)=P(S_n<x\sqrt{B_n})$$
and assume that $S_n$ is well-defined and $B_n<\infty$ for each $n\in\mathbb{N}$.
The following is the main result of this paper.

\begin{theorem} \label{Theorem1}
Suppose that, for all $ n\in\mathbb{N}$ and $j\in\mathbb{Z}^d$, there exist non-negative
constants $c_{nj}$ such that
\begin{equation}  \label{cgf cond}
|L_{nj}(z)|\leq c_{nj}, \quad \forall \, z\in \mathbb C \hbox{  with } |z|<H_n,
\end{equation}
and suppose that $B_n H_n^2\to \infty$ as $n\to \infty$, and
\begin{equation}  \label{cgf cond2}
C_n :=  \sum_{j\in\mathbb{Z}^d}c_{nj} =
O(B_n H_n^2).
\end{equation}
If $x\geq 0 $ and $x=o(H_n\sqrt{B_n})$, then
\begin{equation}  \label{result}
 \frac{1-F_n(x)}{1-\Phi(x)}=\exp\Bigg\{\frac{x^3}{H_n\sqrt{B_n}}
 \lambda_n\Big(\frac{x}{H_n\sqrt{B_n}}\Big)\Bigg\}\Bigg(1+O\Bigg(\frac{x+1}
 {H_n\sqrt{B_n}}\Bigg)\Bigg),
\end{equation}
\begin{equation} \label{result-}
 \frac{F_n(-x)}{\Phi(-x)}=\exp\Bigg\{-\frac{x^3}{H_n\sqrt{B_n}}
 \lambda_n\Big(-\frac{x}{H_n\sqrt{B_n}}\Big)\Bigg\}\Bigg(1+O\Bigg(\frac{x+1}{H_n\sqrt{B_n}}\Bigg)\Bigg),
\end{equation}
where $$\lambda_n(t)=\sum_{k=0}^\infty \beta_{kn}t^k$$
is a power series that stays bounded uniformly in $n$ for sufficiently
small values of $|t|$ and the coefficients $\beta_{kn}$
only depend on the cumulants of $X_{nj}$ $(n \in \mathbb Z, j\in\mathbb{Z}^d)$.
\end{theorem}

For the rest of the paper,  we only state the results for $x\ge 0$. 
Since  $\lambda_n(t)=\sum_{k=0}^\infty \beta_{kn}t^k$ stays bounded uniformly in $n$ for sufficiently small values of $|t|$
and $\beta_{0n}=\frac{H_n}{6B_n}\sum_{j\in\mathbb{Z}^d}\gamma_{3nj}$ from the proof of Theorem \ref{Theorem1},
we have the following corollary:

\begin{corollary}
Assume the conditions of Theorem \ref{Theorem1} hold. Then for  $x\ge 0$ with $x=O\Big((H_n\sqrt{B_n})^{1/3}\Big)$
we have
\begin{align*}
 \frac{1-F_n(x)}{1-\Phi(x)}=\exp\bigg\{\frac{x^3}{6B_n^{3/2}}\sum_{j\in\mathbb{Z}^d}\gamma_{3nj}\bigg\}
 \bigg(1+O\Big(\frac{x+1}{H_n\sqrt{B_n}}\Big)\bigg).
\end{align*}
\end{corollary}
Notice that $\frac{x^3}{6B_n^{3/2}}\sum_{j\in\mathbb{Z}^d}\gamma_{3nj}=O(1)$ under the condition $x=O\Big((H_n\sqrt{B_n})^{1/3}\Big)$.
 Also taking into the account the fact that for $x>0$
$$1-\Phi(x)<\frac{e^{-x^2/2}}{x\sqrt{2\pi}},$$
we obtain the following corollaries:
\begin{corollary} Under the conditions of Theorem \ref{Theorem1}, we have that for  $x\ge 0$ with $x=O\Big((H_n\sqrt{B_n})^{1/3}\Big)$,
\begin{align*}
1-F_n(x)=\Big(1-\Phi(x)\Big)\exp\bigg\{\frac{x^3}{6B_n^{3/2}}\sum_{j\in\mathbb{Z}^d}\gamma_{3nj}\bigg\} +O\bigg(\frac{e^{-x^2/2}}{H_n\sqrt{B_n}}\bigg).
\end{align*}
\end{corollary}
\begin{corollary}
Assume the conditions of Theorem \ref{Theorem1} and $\sum_{j\in\mathbb{Z}^d}\gamma_{3nj}=0$ for all $n\in\mathbb{N}$.
Then for $x \ge 0$ with $x=O\Big((H_n\sqrt{B_n})^{1/3}\Big)$, we have
\begin{align*}
F_n(x)-\Phi(x)=O\bigg(\frac{e^{-x^2/2}}{H_n\sqrt{B_n}}\bigg).
\end{align*}
\end{corollary}
Also as $1-\Phi(x) \sim \frac{1}{x\sqrt{2\pi}}e^{-x^2/2}$, as $x\to\infty$, we have
\begin{corollary} Under the conditions of Theorem \ref{Theorem1}, if $x\to\infty$, $x=o(H_n\sqrt{B_n})$, then
\begin{align*}
\frac{F_n(x+\frac{c}{x})-F_n(x)}{1-F_n(x)}\to1-e^{-c}
\end{align*}
for every positive constant $c$.
\end{corollary}

\section{Applications}\label{application}

In this section, we provide some applications of the main result
in Section \ref{main}. First, we derive a moderate deviation result
for linear random fields with short or long memory; then we apply
this result to risk measures and apply a same argument to study nonparametric regression.

\subsection{Cram\'er type moderate deviation for linear random fields}\label{linearRF}

Let  $X=\{X_j, j\in\mathbb{Z}^d\}$ be a  linear random field defined on
a probability space $(\Omega, {\mathcal F}, P)$ by
$$X_j=\sum_{i\in \mathbb{Z}^d}a_i \varepsilon_{j-i}, \quad j\in\mathbb{Z}^d,$$
where the innovations $\varepsilon_i, i\in\mathbb{Z}^d$, are i.i.d. random variables
with mean zero and finite variances $\sigma^2$, and where $\{a_i, i \in \mathbb{Z}^d\}$
is a sequence of real numbers that satisfy $\sum_{i \in \mathbb{Z}^d} a_i^2 < \infty$.

Linear random fields have been studied extensively in probability and statistics.
We refer to Sang and Xiao (2018) for a brief review on studies in limit theorems,
large and moderate deviations for linear random fields and to Koul et al. (2016),
Lahiri and Robinson (2016) and the reference therein for recent developments in
statistics.

By applying Theorem \ref{Theorem1} in Section \ref{main}, we establish the following
moderate deviation result for linear random fields with short or long memory, under
Cram\'er's condition on the innovations $\varepsilon_i, i\in\mathbb{Z}^d$. Compared
with the moderate deviation results in Sang and Xiao (2018), our Theorem
\ref{Theorem2} below gives more precise convergence rate which holds on much wider
range for $x$.

Suppose that there is a disc centered at $z=0$ within which the cumulant generating
function $L(z) =L_{\varepsilon_i}(z)=\log \EX e^{z\varepsilon_i}$ of $\varepsilon_i$
is analytic and can be expanded in a convergent power series
$$L(z)=\sum_{k=1}^\infty \frac{\gamma_k}{k!}z^k,$$
where $\gamma_k$ is the cumulant of order $k$ of the random variables  $\varepsilon_i, \,
i\in\mathbb{Z}^d$. We have that  $\gamma_1=  \EX\varepsilon_i=0$ and   $\gamma_2
=  \EX\varepsilon_i^2=\sigma^2$, $i\in\mathbb{Z}^d$.

We write
\begin{align}\label{partialsum}
S_n=\sum_{j\in\Gamma^d_n }X_j=\sum_{j\in\mathbb{Z}^d} b_{nj}\varepsilon_j,
\end{align}
where $ b_{nj}=\sum_{i\in\Gamma^d_n}a_{i-j}$. In the setting of Section \ref{main}, we have
$X_{nj} = b_{nj} \varepsilon_j$, $j \in {\mathbb Z}^d$. Then it can be verified that for all $n \ge 1$ 
and $j \in {\mathbb Z}^d$, $X_{nj}$ satisfy condition
(\ref{Cramer}) for suitably chosen $H_n$. In the notation of Section \ref{main}, we have
$$
B_n=\sigma^2\sum_{j\in\mathbb{Z}^d} b_{nj}^2,
\quad\mathrm{  }\quad  F_n(x)=P(S_n<x\sqrt{B_n}).
$$
Hence, we can apply Theorem \ref{Theorem1} to prove the following theorem.
\begin{theorem} \label {Theorem2}
Assume that the linear random field $X= \{X_j, j\in \mathbb{Z}^d\}$ has short memory, i.e.,
\begin{align}\label{Eq:SRD}
A:=\sum_{i\in \mathbb{Z}^d}|a_i|<\infty,  \;\;a:=\sum_{i\in \mathbb{Z}^d}a_i\ne 0,
\end{align}
or long memory with coefficients
\begin{align} \label{Eq:a}
a_i=l(|i|)b(i/|i|)|i|^{-\alpha}, \quad i\in \mathbb{Z}^d,|i|\ne 0,
\end{align}
where $\alpha\in(d/2,d)$ is a constant, $l(\cdot):[1,\infty)\to\mathbb{R}$ is a slowly varying function
at infinity and $b(\cdot)$ is a continuous function defined on the unit sphere ${\mathbb S}_{d-1}$.
Suppose that there exist positive constants $H$ and $C$ such that
\begin{equation}  \label{cgf cond lm}
|L(z)|<C
\end{equation}
in the disc $|z|<H$. Then for all  $x\geq 0$ with $x=o(n^{d/2})$, we have
\begin{align}\label{Eq:lMD}
 \frac{1-F_n(x)}{1-\Phi(x)}=\exp\bigg\{\frac{x^3}{n^{d/2}}\lambda_n\Big(\frac{x}{n^{d/2}}\Big)\bigg\}
 \bigg(1+O\Big(\frac{x+1}{n^{d/2}}\Big)\bigg),
\end{align}
where $$\lambda_n(t)=\sum_{k=0}^\infty \beta_{kn}t^k$$
is a power series that stays bounded uniformly in $n$ for sufficiently small values of $|t|$ and the coefficients $\beta_{kn}$
only depend on the cumulants of $\varepsilon_i$ and on the coefficients $a_i$ of the linear random field.
\end{theorem}

To the best of our knowledge, Theorem \ref{Theorem2} is the first result that gives the
exact tail probability for partial sums of random fields with dependence structure under
the Cram\'er condition.

Due to its preciseness,  Theorem \ref{Theorem2} can be applied to evaluate
the performance of approximation of the distribution of linear random fields
by truncation. We often use the random variable $X_j^m=\sum_{i\in \Gamma_m^d}
a_i \varepsilon_{j-i}$ with finite terms to approximate the
linear random field $X_j=\sum_{i\in \mathbb{Z}^d}a_i \varepsilon_{j-i}$ in
practice. For example, the moving average with finite terms $MA(m)$ is
applied to approximate the linear process (moving average with infinite terms).
In this case, Theorem \ref{Theorem2} also applies to the partial sum $S_n^m
= \sum_{j\in\Gamma^d_n }X_j^m=\sum_{j\in\mathbb{Z}^d} b_{nj}^m\varepsilon_j$.
Here only finite terms $b_{nj}^m$ are non-zero.
Denote
$$
B_n^m=\sigma^2\sum_{j\in\mathbb{Z}^d} (b_{nj}^m)^2,
\quad\mathrm{  }\quad  F_n^m(x)=P(S_n^m<x\sqrt{B_n^m}).
$$
Then for all  $x\geq 0$ with $x=o(n^{d/2})$, we have
\begin{align*}\label{Eq:lMD}
 \frac{1-F_n^m(x)}{1-\Phi(x)}=\exp\bigg\{\frac{x^3}{n^{d/2}}\lambda_n^m
 \Big(\frac{x}{n^{d/2}}\Big)\bigg\}
 \bigg(1+O\Big(\frac{x+1}{n^{d/2}}\Big)\bigg),
\end{align*}
where $$\lambda_n^m(t)=\sum_{k=0}^\infty \beta_{kn}^mt^k,$$
and where the coefficients $\beta_{kn}^m$ have similar definition
as $\beta_{kn}$. To see the difference between the two tail
probabilities of the partial sums, we have
\begin{align*}
 \frac{1-F_n(x)}{1-F_n^m(x)}=\exp\bigg\{\frac{x^3}{n^{d/2}}\Big[\lambda_n
 \Big(\frac{x}{n^{d/2}}\Big)-\lambda_n^m\Big(\frac{x}{n^{d/2}}\Big)\Big]\bigg\}
 \bigg(1+O\Big(\frac{x+1}{n^{d/2}}\Big)\bigg)\\
 =\exp\bigg\{\frac{x^3}{n^{d/2}}\Big[\beta_{0n}-\beta_{0n}^m+\sum_{k=1}^\infty 
 (\beta_{kn}-\beta_{kn}^m)\Big(\frac{x}{n^{d/2}}\Big)^k\Big]\bigg\}
 \bigg(1+O\Big(\frac{x+1}{n^{d/2}}\Big)\bigg),
\end{align*}
here as in the proof of Theorem \ref{Theorem2}, we take $M_n
=\max_{j\in\mathbb{Z}^d}|b_{nj}|$, $H_n=\frac{H}{2M_n}$,
$M_n^m=\max_{j\in\mathbb{Z}^d}|b_{nj}^m|$, $H_n=\frac{H}{2M_n^m}$,
\begin{align*}
\beta_{0n}=\frac{H_n}{6B_n}\sum_{j\in\mathbb{Z}^d}\gamma_{3nj}
=\frac{H\gamma_3}{12M_nB_n}\sum_{j\in\mathbb{Z}^d}(b_{nj})^3,
\end{align*}
\begin{align*}
\beta_{0n}^m=\frac{H_n^m}{6B_n^m}\sum_{j\in\mathbb{Z}^d}\gamma_{3nj}^m
=\frac{H\gamma_3}{12M_n^mB_n^m}\sum_{j\in\mathbb{Z}^d}(b_{nj}^m)^3.
\end{align*}
If $\gamma_3\ne 0$, $\frac{1-F_n(x)}{1-F_n^m(x)}$ is dominated by 
$\exp\big\{\frac{x^3}{n^{d/2}}(\beta_{0n}-\beta_{0n}^m)\big\}$.
If $\gamma_3= 0$, then $\beta_{0n}=\beta_{0n}^m=0$ and
$\frac{1-F_n(x)}{1-F_n^m(x)}$  can be dominated by $\exp\big\{\frac{x^4}
{n^{d}}(\beta_{1n}-\beta_{1n}^m)\big\}$ which depends on whether
$\gamma_4=0$. In general, Theorem \ref{Theorem2} can be applied to
evaluate whether the truncated version $X_j^m$ is a good approximation
to $X_j$ in terms of the ratio  $\frac{1-F_n(x)}{1-F_n^m(x)}$ for $x$
in different ranges which depends on the property of the innovation
$\varepsilon$ and the sequence $\{a_i, i \in \mathbb{Z}^d\}$.

Theorem \ref{Theorem2} can be applied to calculate the tail probability
of the partial sum of some well-known dependent models. For example,
the autoregressive fractionally integrated moving average FARIMA$(p,
\beta,q)$ processes in one dimensional case introduced
by Granger and Joyeux (1980) and Hosking (1981), which is defined as
\begin{equation*}\label{farima}
\phi(B)X_n=\theta(B)(1-B)^{-\beta}\varepsilon_{n}.
\end{equation*}
Here $p, q$ are nonnegative integers, $\phi(z)=1-\phi_1z-\cdots -\phi_pz^p$
is the AR polynomial and $\theta(z)=  1+\theta_1z+\cdots \theta_qz^q$ is the
MA polynomial. Under the conditions that $\phi(z)$ and $\theta(z)$ have no
common zeros, the zeros of $\phi(\cdot)$ lie outside the closed unit disk and
$-1/2<\beta<1/2$, the FARIMA($p,\beta,q$)  process has linear process form
$X_n=\sum_{i=0}^{\infty}{a_i\varepsilon_{n-i}}, \;\;n\in\mathbb{N},$  with
$a_i= \frac{\theta(1)}{\phi(1)}\frac{i^{\beta-1}}{\Gamma(\beta)}+O(i^{-1})$.
Here $\Gamma(\cdot)$ is the gamma function.

\subsection{Approximation of risk measures}

Theorem \ref{Theorem2} can be applied to approximate the risk measures
such as quantiles and tail conditional expectations for the partial
sums $S_n$ in (\ref{partialsum}) of  linear random field
$X=\{X_j, j\in\mathbb{Z}^d\}$. Given the tail probability $\alpha\in(0,1)$,
let $Q_{\alpha,n}$ be the upper $\alpha$-th quantile of $S_{n}$. Namely
$P(S_{n}\geq Q_{\alpha ,n})=\alpha$.
By Theorem \ref{Theorem2}, for all $x\ge 0$ with $x=o(n^{d/2})$,
\begin{align*}\label{Eq:lMD}
 P(S_n>x\sqrt{B_n})=\exp\bigg\{\frac{x^3}{n^{d/2}}\lambda_n
 \Big(\frac{x}{n^{d/2}}\Big)\bigg\}(1-\Phi(x))(1+o(1)).
\end{align*}
We  approximate $Q_{\alpha,n}$ by $x_{\alpha}\sqrt{B_n}$,
where $x=x_{\alpha}=o(n^{d/2})$ can be solved numerically from the equation
\begin{align*}
\exp\bigg\{\frac{x^3}{n^{d/2}}\lambda_n\Big(\frac{x}{n^{d/2}}\Big)\bigg\}(1-\Phi(x))=\alpha.
\end{align*}
The tail conditional expectation is computed as
\begin{align*}
\E(S_{n}|S_{n} &  \geq Q_{\alpha,n})=\frac{Q_{\alpha,n}P
(S_{n}\geq Q_{\alpha,n})+{\int_{Q_{\alpha,n}}^{\infty}P(S_{n}\geq
w)dw}}{{P(S_{n}\geq Q_{\alpha,n})}}\\
&=Q_{\alpha,n}+\frac{\sqrt{B_n}}{\alpha}{\int_{Q_{\alpha,n}/\sqrt{B_n}}^{\infty}\exp\bigg\{\frac{y^3}{n^{d/2}}
\lambda_n\Big(\frac{y}{n^{d/2}}\Big)\bigg\}(1-\Phi(y))dy},
\end{align*}
which can be solved numerically.
The quantile and  tail conditional expectation, which are also called value at risk (VaR) or expected shortfall (ES) in
finance and risk theory, are important measures to model the extremal behavior of  random variables in practice.
The precise moderate deviation results in this article provide a vehicle in the computation of these two measures
of time series or spacial random fields. See Peligrad et al. (2014a) for  a brief review of VaR and ES in the literature
and a study of them when a linear process has $p$-th moment ($p>2$) or has a regularly varying tail with exponent
$t>2$.

\subsection{Nonparametric regression}
Consider the following regression model
$$
Y_{n,j}=g(z_{n,j})+X_{n,j}, \quad j\in \Gamma_n^d,
$$
where $g$ is a bounded continuous function on $\mathbb{R}^m$,
$z_{n, j}$'s are the fixed design points over $\Gamma_n^d \subseteq
\mathbb{Z}^d$ with values in a compact subset of $\mathbb{R}^m$, and
$X_{n,j}= \sum_{i\in\mathbb{Z}^d}a_{i}\varepsilon_{n, j-i}$ is a linear
random field over $\mathbb{Z}^d$, where the i.i.d. innovations $\varepsilon_{n,i}$
satisfy the same conditions as in Subsection \ref{linearRF}. 
The kernel regression estimation for the function $g$ on the basis of sample
pairs $(z_{n,j}, Y_{n,j})$, $j\in \Gamma_n^2 \subset \mathbb{Z}^2$  has been 
studied by Sang  and Xiao (2018) under the condition that the i.i.d. innovations 
$\varepsilon_{n,i}$ satisfy $\|\varepsilon_{n,i}\|_p<\infty$ for some $p>2$ and 
(or) the innovations have regularly varying right tail with index $t>2$. See 
Sang and Xiao (2018) for more references in the literature for
regression models with independent or weakly dependent random field errors.


We study the kernel regression estimation for the function $g$ on the basis 
of sample pairs $(z_{n,j}, Y_{n,j})$, $j\in \Gamma_n^d$, when the i.i.d. innovations 
$\varepsilon_{n,i}$ satisfy the conditions as in Subsection \ref{linearRF}.
Same as in Sang and Xiao (2018) and the other references in the literature, 
the estimator that we consider is given by 
\begin{equation}\label{Eq:K}
g_n(z)=\sum_{j\in\Gamma_n^d} w_{n,j}(z)Y_{n, j},
\end{equation}
where the weight functions $w_{n,j}(\cdot)$'s on $\mathbb R^m$ have form
$$
w_{n,j}(z)=\frac{K(\frac{z-z_{n,j}}{h_n})}
{\sum_{i\in\Gamma_n^d} K(\frac{z-z_{n,i}}{h_n})}.
$$
Here $K: \mathbb{R}^m\rightarrow \mathbb{R}^+$ is a kernel function
and $h_n$ is a sequence of bandwidths which goes to zero as $n
\rightarrow \infty$. Notice that  the weight functions satisfy the 
condition $\sum_{j\in\Gamma_n^d}
w_{n,j}(z)  = 1$.

For a fixed $z \in \R^m$, let 
\begin{equation*}\label{res}
S_n(z) = g_n(z)-\mathbb{E}g_n(z) = \sum_{j\in\Gamma_n^d} w_{n,j}(z)X_{n, j}
=\sum_{j\in\mathbb{Z}^d} b_{n,j}(z)\varepsilon_{n,j},
\end{equation*}
where $b_{n,j}(z)=\sum_{i\in\Gamma_n^d} w_{n,i}(z) a_{i-j}$.
Let $B_n(z)=\sigma^2\sum_{j\in \mathbb{Z}^d}b_{n,j}^2(z)$,  $M_n(z)=\max\limits_{j\in\mathbb{Z}^d}| b_{nj}(z)|$.
 By the same analysis as in the proof of Theorem \ref{Theorem2}, 
we take $H_n\propto  M_n(z)^{-1}$ and derive a moderate deviation result for $S_n(z)=g_n(z)-\mathbb{E}g_n(z)$.  That is,
if $B_n (z)H_n^2\to \infty$ as $n\to \infty$, $x\geq 0, \, x=o(H_n\sqrt{B_n(z)})$, then
\begin{equation}\label{Eq:Tailp}
 P \big(S_n(z)>x\sqrt{B_n(z)}\big)=\left(1-\Phi(x)\right)\exp\bigg\{\frac{x^3}{H_n\sqrt{B_n(z)}}
 \lambda_n\Big(\frac{x}{H_n\sqrt{B_n(z)}}\Big)\bigg\}
 \bigg(1+O\Big(\frac{x+1}{H_n\sqrt{B_n(z)}}\Big)\bigg).
\end{equation}
A similar bound can be derived for $P\big(|S_n(z)|>x\sqrt{B_n(z)} \big)$. 
Notice that these tail probability estimates are more precise than those obtained in Sang and Xiao (2018),
where an upper bound for the law of the iterated logarithm of $g_n(z)-\mathbb{E}g_n(z)$ was derived. With 
the more precise bound on the tail probability in (\ref{Eq:Tailp}) and certain assumptions on $g$ and 
the fixed design points $\{z_{n, j}\}$ [cf. Gu and Tran (2009)], one can construct a confidence interval for 
$g(z)$.   

More interestingly, our method in this paper provides a way for constructing confidence bands for the
function $g(z)$ when $z \in T$, where $T \subset \R^m$ is a compact interval. 
Observe that for any $z, z' \in T$, we can write
\[
S_n(z) - S_n(z') = \sum_{j\in\mathbb{Z}^d} \big(b_{n,j}(z) - b_{n,j}(z')\big) \varepsilon_{n,j}.
\]
Under certain regularity assumption on $g$ and the fixed design points $\{z_{n, j}\}$ [cf. Gu and Tran (2009)],
we can apply the argument in Subsection \ref{linearRF} to derive exponential upper bound for 
the tail probability $P \big(|S_n(z) - S_n(z')|>x\sqrt{B_n(z, z')}\big),$
where $B_n(z, z') = \sigma^2\sum_{j\in \mathbb{Z}^d}\big(b_{n,j}(z) - b_{n,j}(z') \big)^2$. Such a sharp 
upper bound, combined with the chaining argument [cf. Talagrand (2014)] would allow us to derive an 
exponential upper bound for 
$$
P\bigg( \sup_{z, z' \in T} \frac{|S_n(z) - S_n(z')|}{\sqrt{B_n(z, z')}} > x\bigg),
$$
which can be applied to derive uniform convergence rate of $g_n(z) \rightarrow g(z)$ for all $z \in T$ and
to construct confidence band for the function $g(z),  \, z \in T$. It is non-trivial to carry out this 
project rigorously and the verification of the details is a little lengthy. Hence we will have to consider it 
elsewhere.

\section{Proofs}\label{proof}

{\bf Proof of Theorem \ref{Theorem1}}\\

Since $\gamma_{1nj}=0$,  the cumulant generating function $L_{nj}(z)$ of
$X_{nj}$ can be written as
$$
L_{nj}(z)=\log \EX e^{zX_{nj}}=\sum_{k=2}^\infty \frac{\gamma_{knj}}{k!}z^k.$$

Cauchy's inequality for the derivatives of analytic functions together with the condition
\eqref{cgf cond} yields that

\begin{equation}  \label{cumulant}
|\gamma_{knj}|< \frac{k!c_{nj}}{H_n^k}.
\end{equation}

By following the conjugate method (cf. Petrov (1965, 1975)), we now introduce an auxiliary
sequence of independent
random variables $\{\overline{X}_{nj}\}$,  $j\in\mathbb{Z}^d$, with the distribution functions
$$\overline{V}_{nj}(x)=e^{-L_{nj}(z)}\int_{-\infty}^{x}e^{zy}dV_{nj}(y),$$
where $V_{nj}(y)=P(X_{nj}<y)$ and $z \in (-H_n,H_n)$ is a real number whose value will be specified later.

Denote
$$
\overline{m}_{nj}= \EX\overline{X}_{nj}, \quad\mathrm{  }\quad
\overline{\sigma}_{nj}^2= \EX(\overline{X}_{nj}-\overline{m}_{nj})^2,
$$

$$
\overline{S}_{m,n}=\sum_{j\in\Gamma^d_m}\overline{X}_{nj}, \quad\mathrm{  }\quad
\overline{S}_n=\sum_{j\in\mathbb{Z}^d}\overline{X}_{nj},
$$

$$
\overline{M}_{n}=\sum_{j\in\mathbb{Z}^d}\overline{m}_{nj}, \quad\mathrm{  }\quad
\overline{B}_{n}=\sum_{j\in\mathbb{Z}^d}\overline{\sigma}_{nj}^2
 $$
and
$$
\overline{F}_n(x)=P(\overline{S}_n< \overline{M}_{n}+x\sqrt{\overline{B}_{n}}).
$$
Note that, in the above and below, we have suppressed $z$ for simplicity of notations.

We shall see in the later analysis that the quantities $\overline{S}_n, \overline{M}_{n}$ and
$\overline{B}_{n}$  are well-defined for every $n$ and  $z \in \R$ with $|z|<aH_n,$ where $a<1$
 is a positive constant which is independent of $n$. Throughout the proof we will obtain some
 estimates holding for the values of $z$ satisfying $|z|<bH_n,$ where the positive constant
 $b<1$ may vary but is always independent of $n$. We will then take $a$ to be the smallest
 one among those constants $b$. The selection of the constants does not affect the proof
 since the $z=z_n$ we need in the later analysis has property $z=o(H_n)$.

Also, the change of the order of summation of double series presented in the proof is justified
by the absolute convergence of those series in the specified regions. \\

\noindent\textbf{Step 1: Representation of $P(S_{n}<x)$ in terms of the conjugate measure}\\

\noindent First notice that  by equation (2.11) on page 221 of Petrov (1975), for any $m\in \mathbb{N}$, we have
\begin{equation}  \label{partial}
P(S_{m,n}<x)=\exp\bigg\{\sum_{j\in\Gamma^d_m}L_{nj}(z)\bigg\}\int_{-\infty}^{x}e^{-zy}dP(\overline{S}_{m,n}<y).
\end{equation}
Note that the condition \eqref{cgf cond2} implies that  $C_n<\infty, n\in\mathbb{N}$.
From \eqref{cumulant} it follows that for any
$w$ with $|w|<\frac{2}{3}H_n$ and for any  $m\in \mathbb{N}$ we have
\begin{equation}  \label{cgf partial}
\begin{split}
\bigg|\sum_{j\in\Gamma^d_m}L_{nj}(w)\bigg| & = \bigg|\sum_{j\in\Gamma^d_m}\sum_{k=2}^\infty \frac{\gamma_{knj}}{k!}w^k\bigg| \\
&\leq \sum_{j\in\Gamma^d_m}\sum_{k=2}^\infty \frac{|\gamma_{knj}|}{k!}|w|^k\\
&\leq \sum_{j\in\mathbb{Z}^d}\sum_{k=2}^\infty \frac{c_{nj}}{H_n^k}|w|^k \\
&\leq \frac{4}{3}\sum_{j\in\mathbb{Z}^d}c_{nj}=\frac{4}{3}C_n<\infty.
\end{split}
\end{equation}
Therefore,  for any $v$ with $|v|<\frac{1}{2}H_n$ and $z$ with $|z|<\frac{1}{6}H_n$,
\begin{equation}\label{mgf partial}
\begin{split}
&\EX \exp\{v\overline{S}_{m,n}\}=\prod_{j\in\Gamma^d_m} \EX\exp\{v \overline{X}_{nj}\}\\
&=\prod_{j\in\Gamma^d_m}\int_{-\infty}^{\infty}e^{vx}d\overline{V}_{nj}(x) =
\prod_{j\in\Gamma^d_m}\int_{-\infty}^{\infty}e^{vx}e^{-L_{nj}(z)}e^{zx}dV_{nj}(x)\\
&=\prod_{j\in\Gamma^d_m}e^{-L_{nj}(z)}\int_{-\infty}^{\infty}e^{(v+z)x}dV_{nj}(x)= \prod_{j\in\Gamma^d_m}e^{-L_{nj}(z)}e^{L_{nj}(v+z)}\\
&\rightarrow\exp\bigg(\sum_{j\in\mathbb{Z}^d} [L_{nj}(v+z)-L_{nj}(z)]\bigg)<\infty, \;\;\text{as} \;\;m\rightarrow \infty.
\end{split}
\end{equation}
Hence, $\overline{S}_n$ is well-defined and  $\overline{S}_{m,n}$ converges to $\overline{S}_n$ in distribution or equivalently
 in probability or almost surely as $m\rightarrow \infty$.

For the $x$ in $P(S_{n}<x)$, let  $f(y)=\exp\{-zy\}\textbf{1}\{ y<x\}$ and $M>0$. By Markov's inequality,
we have
\begin{align*}
&\EX \Big\{f(\overline{S}_{m,n})\textbf{1}\{ |f(\overline{S}_{m,n})|>M\}\Big\}\\
&\leq \EX \Big\{\exp\{-z\overline{S}_{m,n}\} \textbf{1}\{  \exp\{-z\overline{S}_{m,n}\}>M\}\Big\}\\
&\leq\bigg[\EX \Big\{\exp\{-2z\overline{S}_{m,n}\}\Big\}\bigg]^{\frac{1}{2}}\bigg[\EX \Big\{ \textbf{1}\{ \exp\{-z\overline{S}_{m,n}\}>M\}\Big\}\bigg]^{\frac{1}{2}}\\
&\leq \bigg[ \prod_{j\in\Gamma^d_m}e^{-L_{nj}(z)}e^{L_{nj}(-z)}\bigg]^{\frac{1}{2}}\bigg[\frac{1}{M}\EX \Big\{\exp\{-z\overline{S}_{m,n}\}\Big\}\bigg]^{\frac{1}{2}}\\
&=\frac{1}{\sqrt{M}}\bigg[ \prod_{j\in\Gamma^d_m}e^{-L_{nj}(z)}e^{L_{nj}(-z)}\bigg]^{\frac{1}{2}}\bigg[ \prod_{j\in\Gamma^d_m}e^{-L_{nj}(z)}e^{L_{nj}(0)}\bigg]^{\frac{1}{2}}.
\end{align*}
Hence, by \eqref{cgf partial} we have that for $|z|<\frac{1}{6}H_n$,
$$\lim_{M\to\infty}\limsup_{m\to\infty}\EX \Big\{f(\overline{S}_{m,n})\textbf{1}\{ |f(\overline{S}_{m,n})|>M\}\Big\}=0.$$
Applying Theorem 2.20 from van der Vaart (1998), we have
\begin{equation*}
\int_{-\infty}^{x}e^{-zy}dP(\overline{S}_{m,n}<y)\rightarrow \int_{-\infty}^{x}e^{-zy}dP(\overline{S}_n<y)
\end{equation*}
as $m\rightarrow \infty$. And taking into account that

$$P(S_{m,n}<x)\rightarrow P(S_n<x)$$
and
$$\exp\bigg\{\sum_{j\in\Gamma^d_m}L_{nj}(z)\bigg\}\rightarrow \exp\bigg\{\sum_{j\in\mathbb{Z}^d}L_{nj}(z)\bigg\}$$
as $m\rightarrow \infty$ we obtain from  \eqref{partial} that

\begin{equation}  \label{cdf1}
P(S_n<x)=\exp\bigg\{\sum_{j\in\mathbb{Z}^d}L_{nj}(z)\bigg\}\int_{-\infty}^{x}e^{-zy}dP(\overline{S}_n<y).
\end{equation}

\noindent\textbf{Step 2: Properties of the conjugate measure}\\

\noindent From the calculation of \eqref{mgf partial} it follows that the cumulant generating function $\overline{L}_{nj}(v)$ of the
random variable $ \overline{X}_{nj}$ exists when $|v|$ is sufficiently small and we have
\begin{equation}  \label{Lbar}
\overline{L}_{nj}(v)=-L_{nj}(z)+L_{nj}(v+z),
\end{equation}
$ j\in\mathbb{Z}^d$.
 Denoting by $\overline{\gamma}_{knj}$ the cumulant of order $k$ of the random variable
$ \overline{X}_{nj}$, we obtain

$$\overline{\gamma}_{knj}=\Big[\frac{d^k\overline{L}_{nj}(v)}{dv^k}\Big]_{v=0}=\frac{d^kL_{nj}(z)}{dz^k}.$$

Setting $k=1$ and $k=2$ we find that
\begin{equation} \label{Eq:mbar}
\overline{m}_{nj}=\frac{dL_{nj}(z)}{dz}=\sum_{\ell =2}^\infty \frac{\gamma_{\ell nj}}{(\ell -1)!}z^{\ell-1},
\end{equation}
and
\begin{equation} \label{Eq:sbar}
\overline{\sigma}_{nj}^2=\frac{d^2L_{nj}(z)}{dz^2}=\sum_{\ell=2}^\infty \frac{\gamma_{\ell nj}}{(\ell-2)!}z^{\ell-2}.
\end{equation}
Hence, for $|z|<\frac{1}{2}H_n$, (\ref{Eq:mbar}) imples
\begin{equation} \label{M_nn}
\begin{split}
|\overline{M}_{n}|=\Big|\sum_{j\in\mathbb{Z}^d}\overline{m}_{nj} \Big |=
\Big|\sum_{j\in\mathbb{Z}^d}\sum_{k=2}^\infty \frac{\gamma_{knj}}{(k-1)!}z^{k-1}\Big|\\
\leq\sum_{j\in\mathbb{Z}^d}\sum_{k=2}^\infty \frac{k!c_{nj}}{H_n^k}\frac{|z|^{k-1}}{(k-1)!}\leq\frac{3}{H_n}\sum_{j\in\mathbb{Z}^d}c_{nj}=\frac{3C_n}{H_n},
\end{split}
\end{equation}
which means that $\overline{M}_{n}$ is well-defined and, as a function of $z \in \mathbb C$, is analytic
in $|z|<\frac{1}{2}H_n$.

Also, without loss of generality,
we assume that
\begin{equation} \label{1}
\limsup\limits_n\frac{C_n}{B_n H_n^2}\leq 1.
\end{equation}
By the definition of $\overline{M}_{n}$ and (\ref{Eq:mbar}), we have
\begin{equation} \label{sum for M_n}
\begin{split}
\overline{M}_{n}
&=z\sum_{j\in\mathbb{Z}^d} \gamma_{2nj}+\sum_{j\in\mathbb{Z}^d}\sum_{k=3}^\infty \frac{\gamma_{knj}}{(k-1)!}z^{k-1}\\
&=zB_n+\sum_{j\in\mathbb{Z}^d}\sum_{k=3}^\infty \frac{\gamma_{knj}}{(k-1)!}z^{k-1}.
\end{split}
\end{equation}
It follows from  \eqref{cumulant}  that
\begin{align*}
\begin{split}
\bigg|\sum_{k=3}^\infty \frac{\gamma_{knj}}{(k-1)!}z^{k-1}\bigg| 
&\leq |z|\sum_{k=3}^\infty \frac{k!c_{nj}}{H_n^k}\frac{|z|^{k-2}}{(k-1)!}\\
&=\frac{|z|c_{nj}}{H_n^2}\sum_{k=3}^\infty k\Big| \frac{z}{H_n}\Big|^{k-2}
 \leq \frac{|z|c_{nj}}{2H_n^2}
\end{split}
\end{align*}
for $|z|< b_1H_n$ and a suitable positive constant $b_1<1$ which is independent of $j$ and $n$.
This together with  \eqref{sum for M_n} implies that for $|z|< b_1H_n$
\begin{align*}
|z|\Big(B_n-\frac{C_n}{2H_n^2}\Big) \leq|\overline{M}_{n}|\leq |z|\bigg(B_n+\frac{C_n}{2H_n^2}\bigg).
\end{align*}
Taking into account the condition \eqref{1}, we get that
\begin{equation} \label{speed of M_n}
\overline{M}_{n}\propto |z|B_n.
\end{equation}
Moreover, \eqref{sum for M_n} implies that for $|z|<\frac{1}{2} H_n$,
\begin{equation}   \label{M_n-zB_n}
\begin{split}
\big|\overline{M}_{n}-zB_n \big|
&\leq\sum_{j\in\mathbb{Z}^d}\sum_{k=3}^\infty \frac{k!c_{nj}}{H_n^k}\frac{|z|^{k-1}}{(k-1)!}\\
&\leq\frac{|z|^2}{H_n^3}\sum_{j\in\mathbb{Z}^d}\sum_{k=3}^\infty kc_{nj}\frac{|z|^{k-3}}{H_n^{k-3}}
\leq \frac{8|z|^2C_n}{H_n^3}.
\end{split}
\end{equation}
 Also, by the definition of $\overline{B}_{n}$ and (\ref{Eq:sbar}), we have
\begin{equation}  \label{sum for B_n}
\begin{split}
\overline{B}_{n} 
&=\sum_{j\in\mathbb{Z}^d} \gamma_{2nj}+\sum_{j\in\mathbb{Z}^d}\sum_{k=3}^\infty \frac{\gamma_{knj}}{(k-2)!}z^{k-2}\\
&=B_n+\sum_{j\in\mathbb{Z}^d}\sum_{k=3}^\infty \frac{\gamma_{knj}}{(k-2)!}z^{k-2}.
\end{split}
\end{equation}
It follows from  \eqref{cumulant}   that
\begin{align*}
\begin{split}
\bigg|\sum_{k=3}^\infty \frac{\gamma_{knj}}{(k-2)!}z^{k-2}\bigg|
\leq\sum_{k=3}^\infty \frac{k!c_{nj}}{H_n^k}\frac{|z|^{k-2}}{(k-2)!}\leq\frac{c_{nj}}{2H_n^2}
\end{split}
\end{align*}
for  $|z|< b_2H_n$  and a suitable positive constant $b_2<1$ which is independent of $j$ and $n$.
This together with  \eqref{sum for B_n} implies that for $|z|< b_2H_n$,  $\overline{B}_{n}$ is well-defined and
\begin{align*}
B_n-\frac{C_n}{2H_n^2} \leq|\overline{B}_{n}|\leq B_n+\frac{C_n}{2H_n^2}.
\end{align*}
Condition \eqref{1} then implies that
\begin{equation} \label{speed of B_n}
\overline{B}_n\propto B_n.
\end{equation}

Furthermore, (\ref{sum for B_n}) and  \eqref{cumulant}  imply that for   $|z|<\frac{1}{2}H_n$,
\begin{equation}  \label{B_n tail}
\begin{split}
\Big| \overline{B}_{n}-B_n\Big| 
&\leq\sum_{j\in\mathbb{Z}^d}\sum_{k=3}^\infty \frac{k!c_{nj}}{H_n^k}\frac{|z|^{k-2}}{(k-2)!}\\
&\leq\frac{|z|}{H_n^3}\sum_{j\in\mathbb{Z}^d}\sum_{k=3}^\infty k(k-1)c_{nj}\frac{|z|^{k-3}}{H_n^{k-3}}\\
&\leq \frac{28|z|C_n}{H_n^3}.
\end{split}
\end{equation}

\noindent\textbf{Step 3: Selection of $z$}\\

\noindent Let $z=z_n$ be the real solution of the equation
\begin{equation}  \label{x}
x=\frac{\overline{M}_{n}}{\sqrt{B_n}},
\end{equation}
and let
\begin{equation} \label{t}
t=t_n=\frac{x}{H_n\sqrt{B_n}}.
\end{equation}

Then
\begin{equation}  \label{t1}
\begin{split}
t=\frac{\overline{M}_{n}}{H_nB_n}
=\frac{1}{H_nB_n}\sum_{j\in\mathbb{Z}^d}\sum_{k=2}^\infty \frac{\gamma_{knj}}{(k-1)!}z^{k-1}.
\end{split}
\end{equation}

By \eqref{M_nn} we know that $\frac{\overline{M}_{n}}{H_nB_n}$ is analytic in a disc $|z|<\frac{1}{2}H_n$ and
\begin{align*}
\Big|\frac{\overline{M}_{n}}{H_nB_n}\Big|\leq\frac{3C_n}{H_n^2B_n}
\end{align*}
in that disc.
It follows from Bloch's theorem (see, e.g., Privalov (1984), page 256)  that \eqref{t1} has a real solution which can be written as
\begin{equation}  \label{inverse}
z=\sum\limits_{m=1}^\infty a_{mn}t^m
\end{equation}
for
$$|t|<\bigg(\sqrt{\frac{1}{2}+\frac{3C_n}{H_n^2B_n}}-\sqrt{\frac{3C_n}{H_n^2B_n}} \bigg)^2.$$
Moreover,  the absolute value of that sum in \eqref{inverse} is less than $\frac{1}{2}H_n$.
Condition \eqref{cgf cond2} implies that there exists a disc with center at $t=0$
and radius $R$ that does not depend on $n$ within which the series on the right side of
\eqref{inverse} converges.

It can be checked from \eqref{t1} and \eqref{inverse} that
\begin{equation}\label{a1a2}
a_{1n}=H_n\;\; \text{and}\;\;  a_{2n}=-\frac{H_n^2}{2B_n}\sum_{j\in\mathbb{Z}^d}\gamma_{3nj}.
\end{equation}
Cauchy's inequality implies that for every $m\in\mathbb{N}$,
\begin{align*}
|a_{mn}|\leq\frac{H_n}{2R^m}.
\end{align*}
Therefore,
as $t\to 0$, $a_{1n}t$ becomes the dominant term of the series in \eqref{inverse}.
Hence, for sufficiently large $n$ we have
\begin{align*}
\frac{1}{2}tH_n\leq z\leq 2tH_n,\;\; z=o(H_n)
\end{align*}
and taking into account \eqref{t} we get
\begin{equation}  \label{z}
\frac{x}{2\sqrt{B_n}} \leq z\leq \frac{2x}{\sqrt{B_n}}.
\end{equation}

It follows from \eqref{cgf partial} and \eqref{M_nn} that for $z<\frac{1}{2}H_n$,
$$\Big|z\overline{M}_{n}-\sum_{j\in\mathbb{Z}^d}L_{nj}(z)\Big|\leq \frac{3|z|}{H_n}C_n+\frac{4}{3}C_n<3C_n.$$
For the solution $z$ of the equation \eqref{x} we also have
\begin{equation}  \label{lambda}
\begin{split}
z\overline{M}_{n}-\sum_{j\in\mathbb{Z}^d}L_{nj}(z)&=\sum_{j\in\mathbb{Z}^d}\sum_{k=2}^\infty \frac{\gamma_{knj}}
{(k-1)!}z^{k}-\sum_{j\in\mathbb{Z}^d}\sum_{k=2}^\infty \frac{\gamma_{knj}}{k!}z^k\\
&=\sum_{j\in\mathbb{Z}^d}
\sum_{k=2}^\infty \frac{(k-1)\gamma_{knj}}{k!}\bigg(\sum\limits_{m=1}^\infty a_{mn}t^m \bigg)^k\\
&:=\sum_{j\in\mathbb{Z}^d}\frac{\gamma_{2nj}}{2}a_{1n}^2t^2-\sum_{k=3}^\infty b_{kn}t^k\\
& =\frac{H_n^2B_nt^2} {2}-H_n^2B_nt^3\sum_{k=3}^\infty\frac{ b_{kn}}{H_n^2B_n}t^{k-3}\\
&=\frac{H_n^2B_nt^2}{2}-H_n^2B_nt^3\lambda_n(t),
\end{split}
\end{equation}
where $\lambda_n(t) =\sum_{k=0}^\infty \beta_{kn}t^k$ with $ \beta_{kn} =   b_{(k+3)n}(H_n^2B_n)^{-1}$.

Recall that the series $\sum_{m=1}^\infty a_{mn} t^m$  converges in the disc centered at $t=0$
with radius $R>0$  that does not depend on $n$, and the absolute value of this sum is less than $\frac{1}{2}H_n$.
We see from (\ref{lambda}) that the function $\lambda_n(t)$ is obtained by the substitution of
$\sum_{m=1}^\infty a_{mn}t^m$ in a series that converges on the interval $(-\frac{1}{2}H_n,\frac{1}{2}H_n)$.
It follows from Cauchy's inequality that
$$\big|\beta_{kn} \big| \leq\frac{3C_n}{H_n^2B_nR^{k+3}}\leq\frac{3}{R^{k+3}}, \quad k\geq 0,$$
which means that for $|t|<\frac{1}{2}R$, $\lambda_n(t)$ stays bounded uniformly in $n$.
In particular, by \eqref{a1a2} and \eqref{lambda}, we have  $\beta_{0n}
=\frac{H_n}{6B_n}\sum_{j\in\mathbb{Z}^d}\gamma_{3nj}$.

From now on we will assume that $z$ is the unique real solution of the equation \eqref{x}. \\

\noindent\textbf{Step 4: The case $0\le x\le 1$}\\

\noindent Now we prove the theorem for the case $0\leq x\leq 1$ using the method presented in Petrov and  Robinson (2006).
Throughout the proof, $C$ denotes a positive constant which may vary from line to line, but is independent of $j,n$ and $z$.
If $f_n(s)$ is the characteristic function of $S_n/\sqrt{B_{n}}$ we then have that for $|s|<H_n \sqrt{B_{n}}/2$
\begin{align*}
\begin{split}
&f_n(s)
=\int\limits_{-\infty}^\infty e^{isu}dP(S_n\leq u\sqrt{B_{n}})\\
&=\int\limits_{-\infty}^\infty e^{isy/\sqrt{B_{n}}}dP(S_n\leq y)\\
&=\exp\bigg\{\sum_{j\in\mathbb{Z}^d}L_{nj}(is/\sqrt{B_{n}})\bigg\}.
\end{split}
\end{align*}

Then
\begin{align*}
\begin{split}
&\log f_n(s)= \sum_{j\in\mathbb{Z}^d}L_{nj}(is/\sqrt{B_{n}})=\sum_{j\in\mathbb{Z}^d}\sum_{k=2}^\infty
\frac{\gamma_{knj}}{k!}(is/\sqrt{B_{n}})^k\\
&=-\sum_{j\in\mathbb{Z}^d}\frac{\gamma_{2nj}}{2}s^2/B_{n}+\sum_{j\in\mathbb{Z}^d}\sum_{k=3}^\infty
\frac{\gamma_{knj}}{k!}(is/\sqrt{B_{n}})^k=-s^2/2+\sum_{j\in\mathbb{Z}^d}\sum_{k=3}^\infty \frac{\gamma_{knj}}{k!}(is/\sqrt{B_{n}})^k.
\end{split}
\end{align*}
Thus, using \eqref{cumulant} we get that for $|s|<\delta H_n \sqrt{B_{n}}/2$, with $0<\delta<1$,
\begin{align*}
\begin{split}
|\log f_n(s)+s^2/2|\leq \sum_{j\in\mathbb{Z}^d}\sum_{k=3}^\infty c_{nj}\bigg(\frac{|s|}{H_n\sqrt{B_{n}}}\bigg)^k\leq
C_n\bigg(\frac{|s|}{H_n\sqrt{B_{n}}}\bigg)^3(1-\delta)^{-1}
\end{split}
\end{align*}
Then, for appropriate choice of $\delta$ we have that
\begin{align*}
\begin{split}
|f_n(s)-e^{-s^2/2}|<C\frac{e^{-s^2/4}|s|^3C_n}{H_n^3\sqrt{B_{n}}^3}<C\frac{e^{-s^2/4}|s|^3}{H_n\sqrt{B_{n}}},
\end{split}
\end{align*}
for $|s|<\delta H_n \sqrt{B_{n}}/2$.
Now applying Theorem 5.1 from Petrov (1995) with $b=1/\pi$ and $T=\delta H_n \sqrt{B_{n}}/2$ we get that
\begin{equation} \label{P}
\sup\limits_x|F_n(x)-\Phi(x)|<\frac{C}{H_n\sqrt{B_{n}}}.
\end{equation}

Since $0\le x\le 1$, $B_nH_n^2\rightarrow \infty$ as $n\rightarrow \infty$,  and $\lambda_n\Big(\frac{x}{H_n\sqrt{B_n}}\Big)$
is bounded uniformly in $n$, we have
$$\exp\bigg\{\frac{x^3}{H_n\sqrt{B_n}}\lambda_n\Big(\frac{x}{H_n\sqrt{B_n}}\Big)\bigg\}=1+O(H_n^{-1}B_n^{-1/2}).$$
Together with condition \eqref{cgf cond2}, to have \eqref{result} in the case $0\le x\le 1$, it is sufficient to show
$$ \frac{1-F_n(x)}{1-\Phi(x)}=1+O\bigg(\frac{C}{H_n\sqrt{B_{n}}}\bigg),$$
which is given by \eqref{P}, since $1/2\le \Phi(x)\le \Phi(1)$ for $0\le x\le 1$.

So we will limit the proof of the theorem to the case $x>1, x=o(H_n\sqrt{B_n})$.\\

\noindent\textbf{Step 5: The case $x>1,\; x=o(H_n\sqrt{B_n})$}\\

\noindent Making a change of variables  $ y\rightsquigarrow \overline{M}_n+y\sqrt{\overline{B}_n}$ and applying \eqref{x}, we can rewrite \eqref{cdf1} as

\begin{align}  \label{cdf3}
1-F_n(x)&=\exp\bigg\{-z\overline{M}_n+\sum_{j\in\mathbb{Z}^d}L_{nj}(z)\bigg\}\int_{(x\sqrt{B_n}-\overline{M}_n)/\sqrt{\overline{B}_n}}^{\infty}
\exp\Big\{-zy\sqrt{\overline{B}_n}\Big\}d\overline{F}_n(y)\notag\\
&=\exp\bigg\{-z\overline{M}_n+\sum_{j\in\mathbb{Z}^d}L_{nj}(z)\bigg\}\int_0^{\infty}
\exp\Big\{-zy\sqrt{\overline{B}_n}\Big\}d\overline{F}_n(y).
\end{align}
Denote $r_n(x)=\overline{F}_n(x)-\Phi(x)$ and we show that for sufficiently large $n$
\begin{equation}  \label{r_n}
\sup\limits_x|r_n(x)|\leq \frac{C}{H_n\sqrt{B_{n}}}.
\end{equation}

Let $\overline{f}_n(s)$ be the characteristic function of $(\overline{S}_n-\overline{M}_{n})/\sqrt{\overline{B}_{n}}$.
We then have that
\begin{align*}
\begin{split}
&\overline{f}_n(s)
=\int\limits_{-\infty}^\infty e^{isu}dP(\overline{S}_n\leq u\sqrt{\overline{B}_{n}}+\overline{M}_{n})\\
&=\int\limits_{-\infty}^\infty e^{is(y-\overline{M}_{n})/\sqrt{\overline{B}_{n}}}dP(\overline{S}_n\leq y)\\
&=\exp\bigg\{-is\overline{M}_{n}/\sqrt{\overline{B}_{n}}-\sum_{j\in\mathbb{Z}^d}L_{nj}(z)\bigg\}\int\limits_{-\infty}^\infty e^{(z+is/\sqrt{\overline{B}_{n}})y}dP(S_n\leq y)\\
&=\exp\bigg\{-is\overline{M}_{n}/\sqrt{\overline{B}_{n}}-\sum_{j\in\mathbb{Z}^d}L_{nj}(z)+\sum_{j\in\mathbb{Z}^d}L_{nj}(z+is/\sqrt{\overline{B}_{n}})\bigg\}.
\end{split}
\end{align*}
Then by \eqref{Lbar} for$|z|<\frac{1}{2} H_n$ and $|s|<H_n \sqrt{\overline{B}_{n}}/6$ we have that
\begin{align*}
\begin{split}
&\log\overline{f}_n(s)=-is\overline{M}_{n}/\sqrt{\overline{B}_{n}}+\sum_{j\in\mathbb{Z}^d}\overline{L}_{nj}(is/\sqrt{\overline{B}_{n}})\\
&=-\frac{1}{2}s^2+\frac{1}{6}\big(is/\sqrt{\overline{B}_{n}}\big)^3\Big[\frac{d^3\sum_{j\in\mathbb{Z}^d}\overline{L}_{nj}(y)}{dy^3}\Big]_{y=\theta is/\sqrt{\overline{B}_{n}}},
\end{split}
\end{align*}
where $0\leq |\theta|\leq 1$.
For $|z|<\frac{1}{2} H_n$ and $|s|<\delta H_n \sqrt{\overline{B}_{n}}/6$, with $0<\delta<1$, we have that
\begin{align*}
\begin{split}
&\Big|\Big[\frac{d^3\sum_{j\in\mathbb{Z}^d}\overline{L}_{nj}(y)}{dy^3}\Big]_{y=\theta is/\sqrt{\overline{B}_{n}}}\Big|=\Big|\Big[\frac{d^3}{dy^3}\sum_{j\in\mathbb{Z}^d}\sum_{k=1}^\infty \frac{\overline{\gamma}_{knj}}{k!}y^k\Big]_{y=\theta is/\sqrt{\overline{B}_{n}}}\Big|\\
&=\Big|\sum_{j\in\mathbb{Z}^d}\sum_{k=3}^\infty \frac{\overline{\gamma}_{knj}}{(k-3)!}(\theta is/\sqrt{\overline{B}_{n}})^{k-3}\Big|\\
&\leq \sum_{j\in\mathbb{Z}^d}\sum_{k=3}^\infty k(k-1)(k-2)\frac{c_{nj}}{(H_n/2)^k}\left(s/\sqrt{\overline{B}_{n}}\right)^{k-3}=\frac{48C_n}{H_n^3}\bigg(1-\frac{s/\sqrt{\overline{B}_{n}}}{H_n/2}\bigg)^{-4}\\
&\leq \frac{48C_n}{H_n^3}(1-\delta)^{-4}.
\end{split}
\end{align*}
Thus,
\begin{align*}
\begin{split}
|\log\overline{f}_n(s)+s^2/2|<\frac{8|s|^3C_n}{H_n^3\sqrt{\overline{B}_{n}}^3}(1-\delta)^{-4}.
\end{split}
\end{align*}
Then, for appropriate choice of $\delta$ we have that
\begin{align*}
\begin{split}
|\overline{f}_n(s)-e^{-s^2/2}|<C\frac{e^{-s^2/4}|s|^3C_n}{H_n^3\sqrt{\overline{B}_{n}}^3}<C\frac{e^{-s^2/4}|s|^3}{H_n\sqrt{\overline{B}_{n}}}
\end{split}
\end{align*}
for $|s|<\delta H_n \sqrt{\overline{B}_{n}}/6$.
Now applying \eqref{speed of B_n} and Theorem 5.1 from Petrov (1995)  with $b=1/\pi$ and $T=\delta H_n \sqrt{\overline{B}_{n}}/6$, we have \eqref{r_n}.

By  \eqref{r_n} we have
\begin{equation} \label{int}
\begin{split}
\int_0^{\infty}\exp\Big\{-zy\sqrt{\overline{B}_n}\Big\}\, d\overline{F}_n(y)
&=\frac{1}{\sqrt{2\pi}}\int_0^{\infty}\exp\Big\{-zy\sqrt{\overline{B}_n}-\frac{y^2}{2}\Big\}dy-r_n(0)\\
&\qquad \qquad +z\sqrt{\overline{B}_n}\int_0^{\infty}r_n(y)\exp\Big\{-zy\sqrt{\overline{B}_n}\Big\}dy\\
&=\frac{1}{\sqrt{2\pi}}\int_0^{\infty}\exp\Big\{-zy\sqrt{\overline{B}_n}-\frac{y^2}{2}\Big\}dy+\alpha_n,
\end{split}
\end{equation}
 where $|\alpha_n|\leq \frac{C}{H_n\sqrt{B_{n}}}.$

Denote
$$I_1=\int_0^{\infty}\exp\Big\{-zy\sqrt{\overline{B}_n}-\frac{y^2}{2}\Big\}dy=\psi(z\sqrt{\overline{B}_n})$$
and
$$I_2=\int_0^{\infty}\exp\Big\{-\frac{\overline{M}_{n}}{\sqrt{B_n}}-\frac{y^2}{2}\Big\}dy=\psi(\overline{M}_{n}B^{-\frac{1}{2}}_n),$$
where
$$\psi(x)=\frac{1-\Phi(x)}{\Phi '(x)}=e^{\frac{x^2}{2}}\int_x^{\infty}e^{-\frac{t^2}{2}}dt$$
is the Mills ratio which is known to satisfy
\begin{align*}
\begin{split}
\frac{x}{x^2+1}<\psi(x)<\frac{1}{x},
\end{split}
\end{align*}
for all $x>0$. Hence, by \eqref{z} and \eqref{speed of B_n} we obtain
\begin{align*}
\begin{split}
\frac{\alpha_n}{xI_1}
&=\frac{\alpha_nz\sqrt{\overline{B}_n}}{x}+\frac{\alpha_n}{xz\sqrt{\overline{B}_n}}\\
&\leq C\Bigg(\frac{z\sqrt{B_n}}{xH_n\sqrt{B_n}}+\frac{1}{H_n\sqrt{B_n}xz\sqrt{B_n}}\Bigg)\\
&\leq C\Bigg(\frac{1}{H_n\sqrt{B_n}}+\frac{1}{H_n\sqrt{B_n}x^2}\Bigg)\\
&\leq \frac{C}{H_n\sqrt{B_n}}.
\end{split}
\end{align*}
Hence,
\begin{equation} \label{alpha_n}
\alpha_n=I_1O\Big(\frac{x}{H_n\sqrt{B_n}}\Big).
\end{equation}

For every $y_1<y_2$ we have that $\psi(y_2)-\psi(y_1)=(y_2-y_1)\psi '(u)$, where $y_1<u<y_2$. As for $u>0$,
$|\psi'(u)|<u^{-2}$, then using \eqref{cgf cond2}, \eqref{z}, \eqref{speed of M_n}, \eqref{M_n-zB_n}, \eqref{speed of B_n}
and  \eqref{B_n tail} we get that
\begin{align*}
|I_2-I_1|
&=\Big|\psi '(u)\Big|\Big|\overline{M}_{n}B^{-\frac{1}{2}}_n-z\sqrt{\overline{B}_n}\Big|\\
&\leq \frac{1}{u^2\sqrt{B_n}}\Big|\overline{M}_{n}-z\sqrt{B_n}\sqrt{\overline{B}_n}\Big|\\
&\leq\frac{1}{u^2\sqrt{B_n}}\Big(\Big|\overline{M}_{n}-zB_n\Big|+\Big|zB_n-z\sqrt{B_n}\sqrt{\overline{B}_n}\Big|\Big)\\
&\leq \frac{C}{(\frac{1}{4}x)^2\sqrt{B_n}}\Big(\frac{z^2C_n}{H_n^3}+z\sqrt{B}_n\Big|\sqrt{B_n}-\sqrt{\overline{B}_n}\Big|\Big)\\
&\leq \frac{C}{x^2\sqrt{B_n}}\Big(\frac{x^2C_n}{B_nH_n^3}+\frac{x|B_n-\overline{B}_n|}{\sqrt{B_n}+\sqrt{\overline{B}_n}}\Big)\\
&\leq \frac{C}{x^2\sqrt{B_n}}\Big(\frac{x^2C_n}{B_nH_n^3}+\frac{xzC_n}{H_n^3\sqrt{B_n}}\Big)\\
&\leq \frac{C}{x^2\sqrt{B_n}}\Big(\frac{x^2C_n}{B_nH_n^3}+\frac{x^2C_n}{H_n^3B_n}\Big)=\frac{CC_n}{B_n^\frac{3}{2}H_n^3}\\
&\leq \frac{C}{H_n\sqrt{B_n}}.
\end{align*}
Hence,
$$\frac{|I_2-I_1|}{xI_2}\leq \frac{C}{xH_n\sqrt{B_n}\psi(\overline{M}_{n}B^{-\frac{1}{2}}_n)}=\frac{C}
{xH_n\sqrt{B_n}\psi(x)}<\frac{C}{xH_n\sqrt{B_n}}\frac{x^2+1}{x}<\frac{C}{H_n\sqrt{B_n}},
$$
which means that
\begin{equation} \label{I_2}
\begin{split}
I_1=I_2\Big(1+O\Big(\frac{x}{H_n\sqrt{B_n}}\Big)\Big).
\end{split}
\end{equation}
Finally, combining \eqref{cdf3}, \eqref{x}, \eqref{t}, \eqref{lambda}, \eqref{int} and \eqref{alpha_n} we get
\begin{align*}
1-F_n(x)
&=\exp\bigg\{-\frac{H_n^2B_nt^2}{2}+H_n^2B_nt^3\lambda_n(t)\bigg\}\int_{0}^{\infty}\exp\Big\{-zy\sqrt{\overline{B}_n}\Big\}d\overline{F}_n(y)\\
&=\exp\bigg\{-\frac{x^2}{2}+\frac{x^3}{H_n\sqrt{B_n}}\lambda_n\Big(\frac{x}{H_n\sqrt{B_n}}\Big)\bigg\}\bigg(\frac{1}{\sqrt{2\pi}}I_1+\alpha_n\bigg)\\
&=\exp\bigg\{-\frac{x^2}{2}+\frac{x^3}{H_n\sqrt{B_n}}\lambda_n\Big(\frac{x}{H_n\sqrt{B_n}}\Big)\bigg\}\frac{1}{\sqrt{2\pi}}I_1\bigg(1+O\Big(\frac{x}{H_n\sqrt{B_n}}\Big)\bigg).
\end{align*}
By \eqref{I_2} and the fact that $I_2 = \psi(x)$,  we see that
\begin{align*}
 \frac{1-F_n(x)}{1-\Phi(x)}=\exp\bigg\{\frac{x^3}{H_n\sqrt{B_n}}\lambda_n\Big(\frac{x}{H_n\sqrt{B_n}}\Big)\bigg\}\bigg(1+O\Big(\frac{x}{H_n\sqrt{B_n}}\Big)\bigg).
\end{align*}
This proves \eqref{result}. The proof of (\ref{result-}) follows a same pattern and is omitted.

\medskip

\noindent{\bf Proof of Theorem \ref{Theorem2}}

Since  $\gamma_1=0$, we see that the cumulant generating function $L_{nj}(z)$ of the random variable $ b_{nj}\varepsilon_j,
j\in\mathbb{Z}^d,$  is given by
$$L_{nj}(z)=\log \EX e^{zb_{nj}\varepsilon_j}=\sum_{k=2}^\infty \frac{\gamma_k b_{nj}^k}{k!}z^k.$$
Cauchy's inequality for the derivatives of analytic functions together with the condition \eqref{cgf cond lm} yields that
\begin{equation}  \label{cumulant lm}
|\gamma_k|< \frac{k!C}{H^k}.
\end{equation}
Denote $M_n=\max\limits_{j\in\mathbb{Z}^d}| b_{nj}|$. Then by \eqref{cumulant lm}, for any $H_n$ with  $0<H_n\leq
\frac{H}{2M_n}$
and for any $z$ with $|z|<H_n$ we have

\begin{align*}
\big|L_{nj}(z)\big|  
&\leq\sum_{k=2}^\infty \frac{|\gamma_k|| b_{nj}|^k}{k!}|z|^k
\leq C\sum_{k=2}^\infty \frac{| b_{nj}H_n|^k}{H^k} \\
&= \frac{C}{H} \frac{b_{nj}^2H_n^2}{H-| b_{nj}H_n|}  \leq \frac{2Cb_{nj}^2H_n^2}{H^2}.
\end{align*}
Hence,
$$C_n=\sum_{j\in\mathbb{Z}^d} \frac{2Cb_{nj}^2H_n^2}{H^2}=\frac{2CB_nH_n^2}{\sigma^2H^2}.$$
Then by Theorem \ref{Theorem1}, if $B_n H_n^2\to \infty$ as $n\to \infty$, we have
\begin{equation}  \label{result1}
 \frac{1-F_n(x)}{1-\Phi(x)}=\exp\bigg\{\frac{x^3}{H_n\sqrt{B_n}}\lambda_n\Big(\frac{x}{H_n\sqrt{B_n}}\Big)\bigg\}
 \bigg(1+O\Big(\frac{x+1}{H_n\sqrt{B_n}}\Big)\bigg)
\end{equation}
for $x\geq 0, x=o(H_n\sqrt{B_n})$.

If the linear random field has long memory then we have that (see Surgailis (1982), Theorem 2)
$B_n\propto n^{3d-2\alpha}l^2(n)$. As the function $b(\cdot)$ is bounded, then for $j\in\Gamma^d_n$
we have
\begin{align*}
|b_{nj}|
&\leq C_1\sum_{i\in\Gamma^d_n}l(|i-j|)|i-j|^{-\alpha}\\
&\leq C_1\sum\limits_{k=1}^{2dn}k^{d-1}l(k)k^{-\alpha}
\propto n^{d-\alpha}l(n),
\end{align*}
where we have used the fact (see Bingham et al. (1987) or Seneta (1976)) that for a slowly
varying function $l(x)$ defined on $[1,\infty)$ and for any $\theta>-1$,
$$\int_{1}^{x}y^{\theta}l(y)dy\mathbb{\sim}\frac{x^{\theta+1}l(x)}{\theta+1}, \ \ \ \ \hbox{ as } \ x\rightarrow\infty.$$
It follows from the definition of $a_i$ in  (\ref{Eq:a}) that (for sufficiently large $n$) $M_n=\max\limits_{j\in\mathbb{Z}^d}| b_{nj}|$
is attained at some $j\in \Gamma^d_n$.
Hence, $M_n=O(n^{d-\alpha}l(n))$. We take $H_n\propto n^{-d+\alpha}l^{-1}(n)$  which yields
$$H_n\sqrt{B_n}\propto n^{d/2}.$$
Then the result follows from \eqref{result1}.

If the linear random field has short memory, i.e., $A:=\sum_{i\in \mathbb{Z}^d}|a_i|<\infty,  \;\;a:
=\sum_{i\in \mathbb{Z}^d}a_i\ne 0,$ we can take $M_n=A$ and $H_n=\frac{H}{2A}$.
Moreover, we also have
$$
\sum_{j\in\mathbb{Z}^d}|b_{nj}|\leq \sum_{j\in\mathbb{Z}^d}\sum_{i\in\Gamma^d_n}|a_{i-j}|=(2n+1)^d
\sum_{i\in \mathbb{Z}^d}|a_i|=A(2n+1)^d
$$
and
$$\sum_{j\in\mathbb{Z}^d}|b_{nj}|\geq |\sum_{j\in\mathbb{Z}^d}\sum_{i\in\Gamma^d_n}a_{i-j}|
=(2n+1)^d|\sum_{i\in \mathbb{Z}^d}a_i|=|a|(2n+1)^d,$$
which means that $\sum_{j\in\mathbb{Z}^d}|b_{nj}|\propto n^d$.

As for all $n\in\mathbb{N}$ we have that $|b_{nj}|\leq A$ by the definition of $A$, then
$$\sum_{j\in\mathbb{Z}^d}b_{nj}^2\leq A\sum_{j\in\mathbb{Z}^d}|b_{nj}|\leq A^2(2n+1)^d.$$
On the other hand, for $j\in\Gamma^d_{\left \lfloor{n/2}\right \rfloor }$ we have that $|b_{nj}|>|a|/2$ for sufficiently large $n$.
Hence,
\begin{align*}
&\sum_{j\in\mathbb{Z}^d}b_{nj}^2\geq \sum_{j\in\Gamma^d_{\left \lfloor{n/2}\right \rfloor }}b_{nj}^2\geq \frac{a^2}{4}
\big( 2\left \lfloor{n/2}\right \rfloor+1\big)^d .
\end{align*}
Thus, $\sum_{j\in\mathbb{Z}^d}b_{nj}^2\propto n^d$ and
the result follows from \eqref{result1}.

\end{document}